\documentclass[11pt,letterpaper]{amsart}
\usepackage{amssymb}

\theoremstyle{definition}

\theoremstyle{remark}

\title[Shape matching and Teichm\"uller spaces]{3D shape matching and Teichm\"uller spaces of pointed Riemann surfaces}
\author{Claudio Fontanari and Letizia Pernigotti}
\address{Dipartimento di Matematica \\
Universit\`a degli Studi di Trento  \\
Via Sommarive 14 \\
38123 Povo (Italy)}
\email{fontanar@science.unitn.it, pernigotti@science.unitn.it}

\subjclass[2000]{}

\thanks{}
\begin{document}

\begin{abstract}
Shape matching represents a challenging problem in both information 
engineering and computer science, exhibiting not only a wide spectrum of 
multimedia applications, but also a deep relation with conformal geometry. 
After reviewing the theoretical foundations and the practical issues 
involved in this fashinating subject, we focus on two state-of-the-art 
approaches relying respectively on local features (landmark points) 
and on global properties (conformal parameterizations). 
Finally, we introduce the Teichm\"uller space 
$T_{g,n}$ of $n$-pointed Riemann surfaces of genus $g$ into the realm 
of multimedia, showing that its beautiful geometry provides a natural 
unified framework for three-dimensional shape matching.   
\end{abstract}

\maketitle

As incisively summarized in the survey paper \cite{SM}, shape matching 
deals with transforming a shape and measuring the resemblance with another 
one via some similarity measure, thus providing an essential ingredient 
in shape retrieval and registration. By the way, there seems to be no 
universal and standard definition of what a shape is (see \cite{SM} and 
\cite{Td}): \cite{SM} goes even back to Plato's Meno, where Socrates 
claims in "terms employed in geometrical problems" that "figure is limit 
of solid" (\cite{Plato:380BC}), while \cite{Td} relies on Kendall's 
definition of shape as "all the geometrical information that remains 
when location, scale, and rotational effects (Euclidean transformations) 
are filtered out from an object" (\cite{Ken:77}). 

The authors of \cite{Td} also emphasize the progressive shift of interest 
in the last decade from 2D (where shapes are just silhouettes) to 3D 
(where shapes are embedded surfaces), motivated by spectacular advances 
in the field of three-dimensional computer graphics. Unluckily, as 
remarked in \cite{AS}, most 2D methods do not generalize directly to 
3D model matching, thus inducing a flurry of recent research to focus 
on the specific problem of 3D shape retrieval. It is worth stressing  
that, unlike text documents, 3D models do not allow the simplest form 
of searching by keyword and it is now a rather common belief (see for 
instance \cite{3D}) that an efficient 3D shape retrieval system should 
take into account the principles of the human visual system as disclosed 
by cognitive neurosciences. 

Before turning to the description of a couple of current approaches 
moving along these lines, we wish to mention at least another challenging 
application of three-dimensional shape matching. The so-called computational 
anatomy arises in medical imaging, in particular neuroimaging, and according 
to \cite{LM} it involves comparison of the shape of anatomic structures 
between two individuals, and development of a statistical theory which 
allows shape to be studied across populations. 

A main component in this analysis, after obtaining the individual model 
representations for the subjects being studied, is the establishment of 
correspondence of anatomically homologous substructures between the subjects. 
For example, if we are interested in comparing shape differences between faces 
of two individuals in images, we would like to ensure that the coordinates of 
the left eye in one image correspond to the left eye in the other image. 
On a finer scale, we would like to ensure that the left corner of the 
left eye corresponds appropriately (\cite{LM}). More generally, 
following \cite{PS}, we will define the landmarks of an object as 
the points of interest of the object that have important shape 
attributes. Examples of landmarks are corners, holes, protrusions, 
and high curvature points. 

Landmark-based shape recognition is 
motivated by such a concept of dominant points. It uses landmarks 
as shape features to recognize objects in a scene or to establish 
correspondences between objects, by extending in an optimal way 
over the entire structure the correspondence at a finite subset. 
This last process is called landmark matching (\cite{LM}) and 
originates from the viewpoint of the human visual system, which 
suggests that some dominant points along an object contour are 
rich in information content and are sufficient to characterize 
the shape of the object (see \cite{PS} and the references therein). 
In particular, the paper \cite{LM} presents a methodology and 
algorithm for generating diffeomorphisms of the sphere onto itself, 
given the displacement of a finite set of template landmarks.     

The restriction to a spherical domain is quite natural from the 
point of view of brain imaging, where a basic assumption is   
that the topology of the brain surface is the same as that of a 
crumpled sheet and, in particular, does not have any holes or 
self intersections (\cite{Ot}). It then follows by Riemann
Uniformization (see for instance \cite{Jost:97}, Theorem 4.4.1 (iii))
that there exists a conformal diffeomorphism of such a surface 
of genus zero onto a sphere. More generally, every orientable 
compact embedded surface can be made into a Riemann surface 
with conformal structure (see for instance \cite{Leh:87}, 
Theorem IV.1.1).     

As pointed out in \cite{CI} (see also references therein), the 
visual field is represented in the brain by mappings which are, 
at least approximately, conformal. Thus, to simulate the 
imaging properties of the human visual system conformal image 
mapping is a necessary technique. Moreover, beside its theoretical 
soundness, the application of conformal geometry to the 3D shape 
classification problem presents several practical advantages: 
according to \cite{SC}, the conformal structures are independent 
of triangulation, insensitive to resolution, and robust to noises.  

As far as the genus zero case is concerned, explicit conformal 
flattenings can be obtained by numerically solving the 
Laplace-Beltrami equation (the heart of this procedure consists 
in a finite element reduction to a system of linear equations, 
see \cite{Ot} and \cite{CS}) or by deforming a homeomorphism 
in order to minimize the harmonic energy (both convergence of 
the algorithm and uniqueness of the solution are ensured by 
imposing further constraints, see \cite{CC} and \cite{GZ}).   

A natural generalization involves quasi-conformal mappings, 
which do not distort angles arbitrarily (as it is well-known, 
conformal mappings are angle-preserving): in particular, 
the least-squares conformal mappings (introduced in \cite{LS} 
via a least-squares approximation of the Cauchy-Riemann equations 
and in \cite{GZ} through a discrete version of the harmonic energy 
minimization method) provide a natural solution to 3D 
nonrigid  surface alignment and stitching (\cite{GZ}), at least 
in genus zero.  

In higher genus, a parameterization method based on Riemann 
surface structure has been developed in a series of 
papers by Shing-Tung Yau and collaborators (see in particular 
\cite{BS}). The basic idea is to segment the surface according to 
its conformal structure, parameterize the patching using a 
holomorphic one form, and finally glue them together via 
harmonic maps (\cite{SP}). A comprehensive and self-contained 
account of 3D applications of conformal geometry is now 
available in \cite{ccg}. 

Our (perhaps tendentious) survey of the subject has presented 
three-dimensional shape matching as a classification problem 
for Riemann surfaces carrying (land)marked points. Indeed, following  
\cite{HL:95}, let us fix nonnegative integers $g$ and $n$ 
such that $2g-2+n>0$, a compact connected oriented reference 
surface $S_g$ of genus $g$ and a sequence of $n$ distinct 
points $(x_1, \ldots, x_n)$ on $S_g$. By definition, the 
Teichm\"uller space $T_{g,n}$ is the space of conformal 
structures on $S_g$ up to isotopies that fix $\{x_1, \ldots, x_n\}$
pointwise. It is, in a natural way, a complex manifold of 
dimension $3g-3+n$. 

The idea of classifying 3D shapes according to their conformal 
structure goes back to \cite{SC}, where the general principle 
is stated that in nature it is highly unlikely for different 
shapes to share the same conformal structure. 
Of course by Riemann Uniformization this 
fails in genus $g=0$, but the choice of a suitable number $n>>0$ 
of landmark points allows to address in a uniform way also such 
an exceptional case. Indeed, we believe that the Teichm\"uller space 
$T_{g,n}$ can provide a solid mathematical framework to three-dimensional 
shape matching, by supporting a unified geometric theory of 
landmark matching and conformal parameterizations. 

The case $n=0$ has recently been considered in \cite{J:09} and \cite{fncoord}, 
where the difference between two shapes of the same genus $g$ is measured by 
the Euclidean distance between the corresponding Fenchel-Nielsen coordinates 
in the Teichm\"uller space $T_{g,0}$. In order to compute the Fenchel-Nielsen
coordinates of a shape, the associated surface is conformally deformed along 
the surface Ricci flow, until its Gaussian curvature is $-1$ everywhere. The 
surface is then decomposed into several pairs of hyperbolic pants, i.e., genus 
zero surfaces with three geodesic boundaries. Each pair of hyperbolic pants  
is uniquely described by the lengths of its boundaries and the way of glueing 
different pairs of pants is encoded by the twisting angles between two
adjacent pairs of pants which share a common boundary.

Here we show how to adapt to the case $n > 0$ all steps of the algorithm  
described in \cite{fncoord}, \S 4: 
\begin{enumerate}
\item \label{one}
Compute topological pants decomposition.
\item \label{two}
Compute the hyperbolic metric using Ricci flow.
\item \label{three}
Compute hyperbolic pants decomposition.
\item \label{four}
Compute the Fenchel-Nielsen coordinates.
\end{enumerate}  

In order to approximate a punctured Riemann surface by a triangular mesh, we consider 
the corresponding compact surface with boundaries obtained by removing open discs of 
sufficiently small but strictly positive radius instead of points.

Step (\ref{one}) follows the approach of \cite{constmap}, which relies on \cite{handletunnel} 
and is explicitly presented in the case of surfaces with boundary. 

As far step (\ref{two}) is concerned, we recall that the exponential convergence of the solution 
of the Ricci flow equation to a complete metric with prescribed constant scalar curvature has 
been proved in \cite{Chow} for compact surfaces with negative Euler characteristic. The case 
of interest for us, namely, an open surface conformal to a punctured compact Riemann surface, 
has been recently addressed in \cite{riccicusp} (see in particular Theorem 3). 
On the other hand, the combinatorial 
version of the Ricci flow is introduced in \cite{combricci} for triangulations on compact
surfaces with boundary (see in particular Theorem 5.1). Since the description of step (\ref{two})
in \cite{fncoord}, \S 4, relies on \cite{combricci}, it directly extends to our setting 
(for implementation issues in the case of surfaces with boundary we refer to \cite{ccg}, 
\S 12). 
 
Finally, steps (\ref{three}) and (\ref{four}) in \cite{fncoord}, \S 4, rely on 
\cite{discrete}, which addresses also surfaces with boundary, and can be repeated 
verbatim to produce (up to renumbering) a sequence of lenghts and twisting angles 
$$
\{(l_1, \theta_1), (l_2, \theta_2), \ldots , (l_{3g-3+n}, \theta_{3g-3+n}), 
l_{3g-3+n+1}, \ldots, l_{3g-3+2n} \}.
$$
Indeed, there are $3g-3+2n$ length parameters, one for each curve of the pants decomposition
and one for each boundary curve, and $3g-3+n$ twist parameters, one for each curve of the pants 
decomposition. By setting the length parameters $l_{3g-3+n+1}, \ldots, l_{3g-3+2n}$ to be zero, 
we can turn boundary components into punctures (see for instance \cite{primer}, \S 10.6.3). 
The corresponding Fenchel-Nielsen coordinates are given by 
$$
\{(l_1, \theta_1), (l_2, \theta_2), \ldots , (l_{3g-3+n}, \theta_{3g-3+n})\}.
$$

Full implementation details and numerical results will be presented elsewhere. 

All the above undoubtedly points towards a potentially fruitful 
interaction between two such apparently unrelated fields as moduli 
spaces and information technologies. We hope to see deeper into 
this in the next future.

\end{document}